\documentclass[12pt,reqno]{amsproc}

\usepackage{bbm}
\usepackage{amssymb,amsmath}
\usepackage{epsfig}

\usepackage{epsfig}
\newcommand{\ind}{\displaystyle\mathbbm{I}}
\usepackage{amsmath,amsthm,amssymb,bm}
\usepackage{graphicx, color,blindtext}

\usepackage{subfig}
\usepackage{float}

\newcommand{\RR}{\mathbb{R}}

\newcommand{\NN}{\mathbb{N}}

\newcommand{\D}{\ensuremath{\mathcal{D}}}

\newcommand{\T}{\mathcal T}
\newcommand{\C}{\mathcal C}

\newtheorem{thm}{Theorem}[section]
\newtheorem{lem}[thm]{Lemma}

\newtheorem{pro}[thm]{Proposition}
\newtheorem{rem}[thm]{Remark}
\newtheorem{defi}[thm]{Definition}

\everymath={\displaystyle}

\allowdisplaybreaks

\newcommand{\cs}{$^\dagger$} \newcommand{\cm}{$^\ddagger$}

\title[On a  Mathematical model for traveling sand dune]
{On a  Mathematical model for traveling\\  sand dune}

\author[Igbida]{Noureddine Igbida\cs} 
\thanks{\cs Institut de recherche XLIM-DMI, UMR-CNRS 6172, Facult\'e des Sciences et Techniques, Universit\'e de Limoges, France.
E-mail: {\tt noureddine.igbida@unilim.fr} }

\author[Karami]{Fahd Karami\cm} 
\thanks{\cm Ecole Sup\'erieure de Technologie d'Essaouira, Universit\'e Cadi Ayyad, B.P. 383 Essaouira El Jadida, Essaouira, Morocco
E-mail: {\tt fa.karami@uca.ma}, {\tt driss.meskine@laposte.net} }

\author[Meskine]{Driss Meskine \cm}

\date{\today}
\begin{document}

\begin{abstract}
Our aim in this  paper  is to introduce and study a mathematical model for the description of traveling sand dunes.  We use  surface flow process of  sand   under the effect of  wind and gravity. We model this phenomena  by a non linear  diffusion-transport  equation coupling the effect of transportation of sand due to the wind and the avalanches due to the gravity and the repose angle. The avalanche flow  is governed by the evolution surface model and we use a nonlocal term to handle the transport of sand face to the wind. 

 \end{abstract}



\maketitle

%%%%%%%%%%%%%%%%%%%%%%%%%%%%%%%%%%%%%%%%%%%%%%%%%%%%%

\begin{center}
 \section{Introduction}
\end{center}
%%%%%%

    The diffusion-transport  equation of the type 
    \begin{equation}\label{geneq} 
     \partial _t u = \nabla \cdot \Big( m\: \nabla u - V\:  u  \Big)+ f 
     \end{equation}
    governs the spatio-temporal dynamics of a density $u(x,t) $  of particles.  The first term on the right hand side  describes random
motion and the parameter $m$, corresponding to diffusion coefficient, is connected to rate of random exchanges between 
    the particles at the position $x$ and neighbors positions. The second
term is a transportation with velocity  $V,$  and is connected to the transport of particles according  to the 
    the vector field $V.$ 
    
  \medskip

     Our aim here is to show how one can use this type of equation to describe the movement of traveling sand dunes ; the so called Barchans.
  A Barchan is a dune of the shape of a crescent lying in the direction of the wind. It arises where the supply of sand is low and under unidirectional winds.   
 The wind rolled the sand back to the slope of the dune back up the ridge and comes to  cause small avalanches on steep slopes over the front. This furthered the dune.  Such dunes  (Barchan) move in the desert at speeds depending on their size and strength of wind. In one year, they argue a few meters to a few tens of meters.  A wind of $\displaystyle 25 km/h$ is  enough to  turn on the preceding process and furthered the dune.  Our main  is to show how one can couple  between the action of the wind and the action of gravity into the equation (\ref{geneq}) to fashion a model for traveling sand dunes like Barchans. 
  \medskip

Recall that the surface evolution  model    is 
a very useful model for the description of the dynamic of a granular matter under the gravity effect. In \cite{Pri1} (see also \cite{ArEvWu},  and \cite{DuIg1}) the authors show that this toy model  gives a simple way to study the dynamic of the sandpile from the theoretical and a numerical point of view.  In particular, the connection of this model with a stochastic model for sandpile (cf. \cite{EvRe}) shows how it is able to handle the global dynamic of a structure of granular matter structure using simply the repose angle.  In this paper, we  show how one can use  the surface evolution model in the equation (\ref{geneq}) to 
 describe the dynamic of traveling sand dunes.  There is a wide literature concerning mathematical and physical studies of  dunes  (cf.  \cite{Bagnold}, \cite{KoLa1}, \cite{AnClDo1,AnClDo2}, \cite{SaKrHe} and  \cite{KrSaHe}). We do believe that   this is an intricate  phenomena with many 
determinant parameters.    We are  certainly neglecting  some of them.  Nevertheless,  as for the case of  sandpile, we do believe that  simple models (like  toys model) may help to encode  complex phenomenon related to the granular matter structures.  
 
   \medskip   

    In the following section,  we give some preliminaries and present our model for a traveling sand dune. Section 3 is devoted main results   of   existence and uniqueness of a weak solution.

 \bigskip

 %%%%%%%%%%%
 \section{Preliminaries and modeling} 
 \bigskip
 %%%%%%%%%%%%%%%%%%%%%%

 \begin{figure}[!h]  
 \includegraphics[width=12cm,height=8cm]{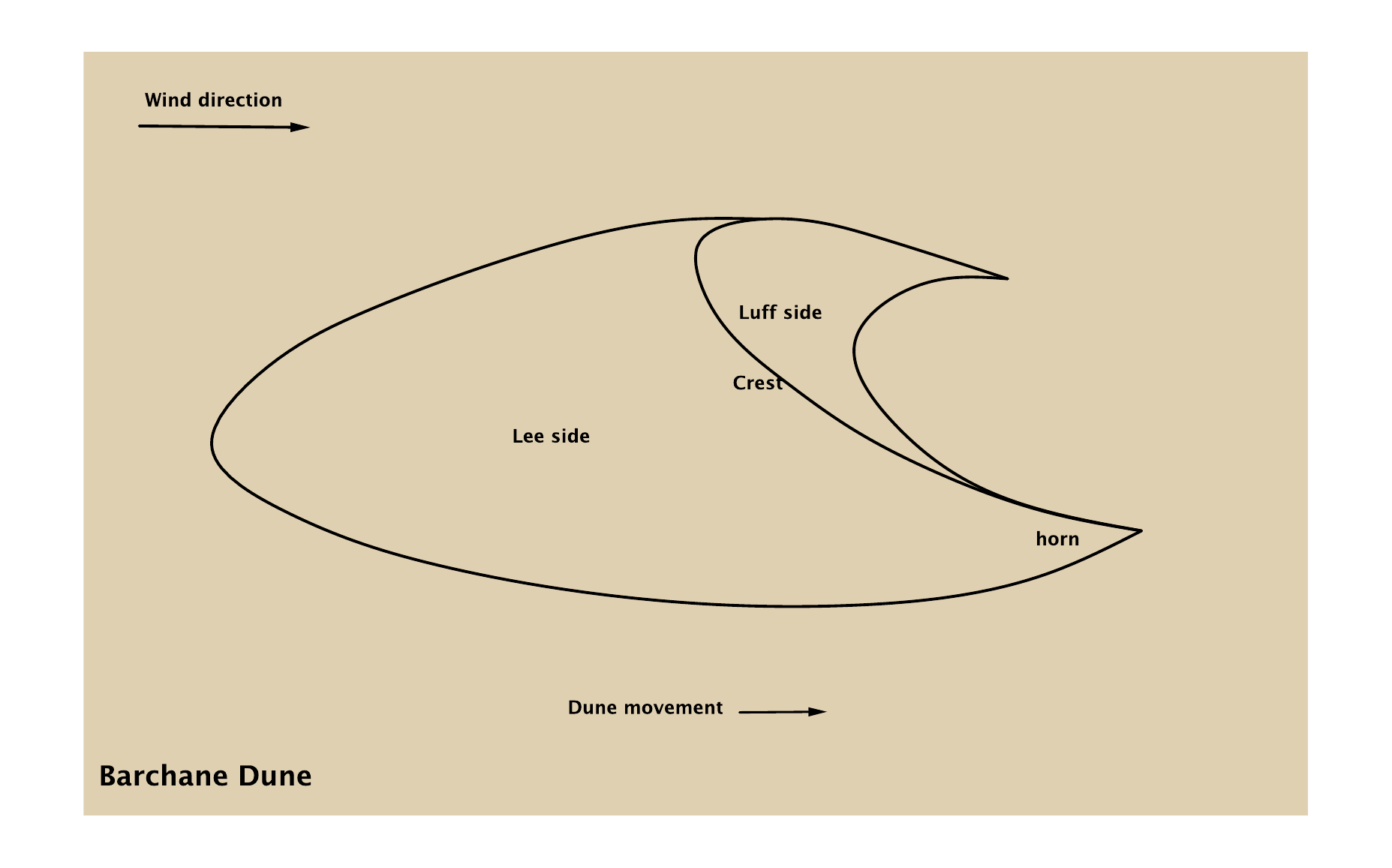}
\caption{\label{Barchan} Barchan}
\end{figure} 

The simplest and the most well-known type of dunes is the Barchan. In general, it has the form of a hill with a Luff side and a Lee  side separated by a crest (cf. Figure \ref{Barchan}). On the Lee side, sand is taken up by the wind into a moving layer,  transported up to the crest and pass to the other side; the Luff  side. By neglecting the effects of precipitation and swirls, on the Luff side, the dynamic of the sand   is  generated by the action of    gravity on the sand arriving at the crest.

 Though there are many speculations and experimental observations on the evolution of the shape, the height and  the distribution of dunes, there is no universal model for the study     of the motion  of  sand dunes.   There is a large literature on this subject  (cf.  \cite{Bagnold}, \cite{KoLa1}, \cite{AnClDo1,AnClDo2}, \cite{SaKrHe} and  \cite{KrSaHe}).  Some mathematical  models treats the  dunes as an aerodynamic objects with an adequate  smooth shape  to let the
air flow by with the least effort (cf.  \cite{KoLa1}).   The so called BCRE models (cf.  \cite{BCRE}) use  the conservation of  mass and the repose angle to build a system of two coupled differential equations  for the height of the topography $h$ and the amount of mobile particles $R.$ The particles are supposed to
move all with the same velocity $u.$   Other simplified and realistic physical models   (cf.  \cite{AnClDo1,AnClDo2}, \cite{SaKrHe} and \cite{KrSaHe})  use the mass and momentum conservation in presence of erosion and external  
 forces  
 %average density and velocity of sand. 
to  derive  coupled differential equations  to study the evolution of  the morphology of dunes. Among other things, these models  focus on  the way in which a sand movement could be constructed from wind data (the choice of formula for linking wind velocity to sand movement, the choice of a threshold velocity for sand movement e.g. Bagnold, etc). Keeping in mind that the driving force for the Barchans is the wind, our aim here is to introduce and study a simple model which combines the effort of wind on the Lee side with the avalanches generated by the repose angle.     More precisely, we introduce and study  a new mathematical model  in the form of diffusion-transport equation   (\ref{geneq})  for the evolution of a morphology of a Barchan under the effect of a unidirectional wind.   

 Let us denote by $u=u(t,x,y)$   the height  of  the dune at time $t\geq 0$  and  at the position $\displaystyle (x,y) $   in the plane $\RR^N$ ($N=2$ in practice). Then, $u$ can be described by (\ref{geneq}) where the term $-m\: \nabla u$   is connected   to the net flux of   the avalanches of sands resulting from the action of the gravity  and the repose angle.   The term   $uV$ is connected to the transport of the sand under the action of wind up to  the crest.    See here that $f\equiv 0$, since we are assuming that  there is no source of sand. 
%   it appears that a quantitative framework for understanding
%the natural behavior concerning the coupling between the transport phenomena and the gravity for modeling  traveling sand dune is still lacking. In general   it is not immediately obvious how to quantify and couples them.  
%
Thanks to \cite{Pri1} (see also \cite{ArEvWu} and \cite{EvRe}) we know that the avalanche can be governed by   a non standard diffusion parameter $m$ (unknown) that is connected to the sub-gardient constraint in the following way 
 $$   m\geq 0,\ \vert \nabla u\vert \leq  \lambda,\    m\: (\vert \nabla u\vert - \lambda )=0, $$
 where  $ \theta =\arctan (\lambda)$  is  the repose angle  of the sand. This is  the consequence of the fact that   the inertia is neglected, the surface flow is directed towards the
steepest descent, the surface slope of the sandpile cannot exceed the repose angle $\theta$ of the material and 
 there is no pouring over the parts of angle less than $\theta.$
 As to the action of the wind,  the resulting phenomena is a transportation   of $u$  in the direction of  the  wind.   Indeed, the wind
proceed by taking up the sand into a  moving layer and transport it  up to  the crest.  This creates a ripping curent of sand concentrated in the Lee side ; face up the wind. 
To handle  the way a sand movement could be constructed from wind data, we use  the nonlocal interactions between the positions of the dune face to the wind.  We assume  that   the velocity $V=(V_1,0)$ of the transported layer  is induced at a site $x$ by the net effects of the slope  of all particles at various sites $y$ around $x.$  More precisely, we  consider  $V_1$ in  the  form 
 $$ V_1=    \mathcal{H}  (K\star \partial_x u) =  \mathcal{H}  \left( \int_{B(x,r)}  K(x-y)\: \partial_x  u (y,t)\: dy\right) , $$
 where  the kernel $K$ associates
a strength of interaction per unit density with the distance $x-y$ 
 between any two sites  over some finite
domain $B(x,r).$    Taking the average of this form may give
more weight to information about particles  that are closer, or those
that are farther away.  Moreover, since the transport need to be restricted to the region  facing the wind, we assume that  the function, $\mathcal{H} \: :\: \RR\to \RR^+$ is a phenomenological parameter which vanishes in   the  region $(-\infty,0).$ 
%So, the net flux corresponding to the wind  transport of sand may be describe by the following  vector field   
%\begin{equation}\label{fluxw}
% (     u\:  \mathcal{H}  (K\star \partial_x u) , 0)
%\end{equation}
% where   $K$ is a given Kernel compactly supported in $B(0,r),$ for a given parameter $r>0.$   
 
\bigskip
Thus, we consider the following model to describe     the evolution of the morphology of a Barchan under the effect of a unidirectional wind  (in the direction $ (1,0)$)  :   
\begin{equation}\label{nlmodel}
\left\{  \begin{array}{ll}   
\left.  \begin{array}{l}  
 \displaystyle  \partial_t u  -\nabla\cdot   ({  m} \:\nabla u) +    \partial_x \Big( u \: \mathcal{H}  (K\star \partial_x u) \Big)   =0 
  \\ \\  
 \displaystyle \vert  \nabla u\vert   \leq   \lambda,   \   {  m}\geq 0,\    
{ m}\: (\vert \nabla u\vert -\lambda)=0
     \end{array} \right\} & \quad\hbox{ in }(0,\infty)\times\RR^N\\  \\  
 u(0,x)=u_{0}(x) &\quad\hbox{ for  }x\in \RR^N.    \end{array}\right.  
\end{equation} 
% where  $u_{0}$  is given (the initial profile of the dune), $\displaystyle \gamma \: :\: \RR^+\to \RR^+$ is  a   continuous function and  $\mathcal{H} \: :\: \RR\to \RR^+$ is compactly supported in $(0,\infty).$  

 \begin{rem} 
 \begin{enumerate}
  
 \item In (\ref{nlmodel}),  we are assuming  that the flux du to the wind  ; i.e. the quantity of sand transported by unit of time through fixed vertical line, depends on  the speed of the wind  and the angle of the position. Moreover, thanks to the assumption on $\mathcal H$ (vanishing in $(-\infty,0)$),   the action of the wind is null whenever the slope is not face to the wind.   To be more general, it is possible to assume that this flux depends also  on the eight, i.e. we can assume that 
 $$V_1=   \: \beta(u) \:  \mathcal{H}  (K\star \partial_x u),$$
 where  $\displaystyle \beta \: :\: \RR^+\to \RR^+$ is  a   continuous function. 
 % and  $\mathcal{H} \: :\: \RR\to \RR^+$ is compactly supported in $(0,\infty).$  

   \item We see in the formal model  (\ref{nlmodel}) that $\mathcal H$ is null for negative values and strictly positive on $\RR^+.$ So, formally,  $\mathcal H=\chi_{\RR^+}.$  Here, for technical reason, we consider continuous approximation of this kind of profile by assuming that  $\mathcal H$ is a    Lipschitz continuous function on $\RR$. 
For instance, one can take   $$ \displaystyle  \mathcal{H} (r) =1-\frac{1}{\sqrt{\pi}}  \int^{-r/\sqrt{\varepsilon}}_{-1/\varepsilon} e^{-z^2}  dz , $$
 where  $0<\varepsilon<1$   is a given fixed parameter $.$

  \item Since the assumptions on $\mathcal H,$ one sees that the crest, which corresponds here to the region where $u$ changes its monotonicity, constitutes a free boundary separating the region of avalanches and the region of wind erosion of sand. Indeed, the transport term $uV$ disappears in the region where $u$ is nonincreasing.

   \item   It is possible to improve the property of  $\mathcal H$ to better  describe the movement of sand face the wind. For instance if we assume that the grains move more and more slowly whenever they are face to important slope, then    $\mathcal H $ can be assumed to be a nondecreasing Lipschitz continuous function.  Typical example may be given by 
 $$  \mathcal H (r)=\frac{r^+}{ \sqrt{1+r^2}},\quad \hbox{ for any } r\in \RR.  $$
In this paper, we just assume that $\mathcal H $ is a Lipschitz continuous function. 
The discussions concerning concrete assumptions on  $\mathcal H $ and also on $\gamma$ and $K$ will be discussed in forthcoming papers.

\item Replacing  $K$ by  $K_\sigma$  in  (\ref{nlmodel}) where $K_\sigma \in \D(\RR^N)$  is a smoothing kernel satisfying  
 \begin{itemize} 
 \item  $\displaystyle \int_{\RR^N} K_\sigma  (x) \:  dx=1$ 
 \item $\displaystyle K_\sigma (x) \to \delta_x$  as $\sigma \to 0$, $\delta_x$ is the Dirac function at the point $x$,
 \end{itemize} 
and letting formally $\sigma \to 0 $ in  (\ref{nlmodel}), we obtain  the following PDE : 
 \begin{equation}\label{pdemodel}
\left\{  \begin{array}{ll}   
\left.  \begin{array}{l}  
 \displaystyle  \partial_t u  -\nabla\cdot   ({  m} \:\nabla u) +   \partial_x \Big( \gamma (u) \: \mathcal{H}  ( \partial_x u) \Big)   =0 
  \\ \\  
 \displaystyle \vert  \nabla u\vert   \leq   \lambda,\ \exists\   {  m}\geq 0,\    
{ m}\: (\vert \nabla u\vert -\lambda)=0
     \end{array} \right\} & \quad\hbox{ in }(0,\infty)\times\RR^N\\  \\  
 u(0,x)=u_{0}(x) &\quad\hbox{ for  }x\in \RR^N.   \end{array}\right.  
\end{equation} 
Since  $\mathcal H $ is assumed to be nondecreasing, one sees that the term $  \partial_x \Big( \gamma (u) \: \mathcal{H}  ( \partial_x u) \Big)$  is a anti-diffusive and may creates some obstruction to the existence of a solution.  It is not clear for us if (\ref{pdemodel})  is well posed in this case or not. 

\end{enumerate}

 \end{rem}

 %%%%%%%%%%%
 \section{Existence and uniqueness} 

 %%%%%%%%%%%%%%%%%%%%%%

  To study the model (\ref{nlmodel}), we restrict our-self to $\Omega\subset \mathbb{R}^N$  a bounded open domain, with Lipschitz boundary $\partial \Omega$ and outer unit normal $\eta$.  Consider the following  nonlocal equation  with Dirichlet boundary condition   :
$$\eqno (E) \label{nlmodelb} 
\left\{  \begin{array}{ll}   
\left.  \begin{array}{l}  
 \displaystyle  \partial_t u  -\nabla\cdot   ({  m} \:\nabla u) +   \partial_x \Big( \gamma (u) \: \mathcal{H}  (K\star \partial_x u) \Big)   =f
  \\ \\  
 \displaystyle \vert  \nabla u\vert   \leq   \lambda,\ \exists\   {  m}\geq 0,\    
{ m}\: (\vert \nabla u\vert -\lambda)=0
     \end{array} \right\} &  \mbox{ in}\,\Omega_T:=(0,T)\times \Omega \\ \\
%\Theta \geq 0,  \; \;\; \Theta (  |\nabla u |- \lambda)=0 \quad \mbox{in}  
 u=0  & \mbox{ in}\,\,\Sigma_T:=(0,T)\times \partial \Omega\, \\\\
u(0,x)=u_0(x)  & \mbox{ in}\,\, \Omega,
\end{array} \right. $$
where $u_0$ patterns  the intial shape of the dune. Here and  throughout the paper, we assume that 

\begin{itemize} 
\item $\mathcal{H} \: :\: \RR\to \RR^+$  is a    Lipschitz continuous function.
\item  $K$ is a given regular Kernel compactly supported in $B(0,r),$ for a given parameter $r>0.$ 
\item $\displaystyle \gamma \: :\: \RR^+\to \RR^+$ is  a    Lipschitz continuous  function with  $\displaystyle \gamma (0)=0$. 
 
\end{itemize}

\bigskip

Our main results concerns existence and uniqueness of a solution.  As usual for the first differential operator governing the PDE (E), we use the notion of variational solution.   We consider $\C_0(\Omega),$ the set of continuous function null on the boundary. For any $0\leq \alpha\leq 1,$ we consider  
$$ \displaystyle \C_0^{0,\alpha}(\Omega)=\Big\{ u\in \C_0(\Omega)\: :\:  u(x)-u(y)\leq C \vert x-y\vert^\alpha \hbox{ for any }x,y\in \overline \Omega   \Big\} ,$$
endowed with the natural norm 
$$\Vert u\Vert _{\C_0^{0,\alpha}(\Omega)} = \sup_{x\in \Omega} \vert u(x)\vert + \sup_{x\neq y\in \Omega} \frac{\vert u(x)-u(y)\vert}{\vert x-y\vert^\alpha} .  $$
We denote by 
$$Lip= \C^{0,1}(\Omega)\quad \hbox{  and }Lip_0= Lip\cap \C_0(\Omega) .$$ 
Then, we denote by  
 $$ \displaystyle Lip_1=\Big\{ u\in \C_0(\Omega)\: :\:  u(x)-u(y)\leq \vert x-y\vert \hbox{ for any }x,y\in   \Omega   \Big\} .$$   
  The topological dual space of $Lip_0$ will be denoted by $Lip^*_0$ and is endowed with the natural dual norm, and we denote by $\langle .,.\rangle$ the duality bracket.  It is clear that, for any $\xi\in Lip_1,$ we have 
  $$  \vert \xi(x)\vert \leq \delta_\Omega, \quad \hbox{ for any }x\in \Omega,$$
  where $\delta_\Omega$ denotes the diameter of the domain $\Omega.$ 
 
 \bigskip 
Recall that the notion of solution for the problems of the type (\ref{nlmodelb}) is not standard in general.   The problem presents two specific difficulties. The first one is  related  to the main operator governing the equation : $-\nabla \cdot (m\nabla u)$ with $m\geq 0$ and  $m(\vert \nabla u\vert -1)=0$. And the second one is connected to the   regularity of the term $\partial _t u$.

Concerning  the main operator governing the equation recall that $u$ is Lipschitz and,  in general even in the  case where  $\mathcal{H} \equiv 0,$ $m  $ is singular.     So, the term $m\nabla u$ is not well  defined in general and needs to be specified.  To handle the PDE  with the operator in divergence form  request the use of the notion of tangential gradient with respect to a measure (cf. \cite{BBS,BBS97,BCJ}). 
Nevertheless, to avoid all the   technicality related to this approach, we use here the notion of truncated-variational solution (that we call simply variational solution) to handle the problem. Indeed, the following lemma strips the way to this alternative.  For any $k>0,$ the real function $T_k$ denotes the usual truncation given by 
\begin{equation}\label{trunc}
\displaystyle T_k(r)=\max(\min(r,k),-k),\quad \hbox{ for any }r\in \RR.   
\end{equation} 
 
% Lemma
 \begin{lem}\label{ldefsol}
 	Let $\eta \in Lip^*_0$ and $u\in Lip_1.$  If, there exists $m\in L^1(\Omega)$ such that $m\geq 0$,  $m(\vert \nabla u\vert -1)=0$ a.e. in $\Omega$ and $\nabla \cdot (m\nabla u)=\eta $ in $\mathcal D'(\Omega),$ then,	
 	$$\langle \eta ,T_k(u-\xi) \rangle \geq 0, \quad \hbox{   for any  }  \xi\in Lip_1 \hbox{ and }k> 0.  $$   
 	\end{lem}
 The proof of this lemma is  simple, we let it as an exercise for the interested reader. Let us notice that the converse part remains true if one take on  $m$ to be a measure  and the gradient to  be the tangential gradient with respect to $m.$ Other equivalent formulations may be found in \cite{Igmonge1}.

This being said, one sees that performing the notion of variational solution  in (E)  generates formally the quantity $ \langle u_t  ,T_k(u-\xi)\rangle.$   Since in general  $\partial_t u$  is not necessary a Lebesgue function we  process the following (formal) integration by parts formula in the definition of the solution :  
$$\langle  \partial_ t u, T_k(u-\xi) \rangle  =
\frac{d}{dt}  \int_{\Omega}  \int_0^{u(t)} T_k (s -\xi) ds  dx.$$
Observe that,  letting $k\to\infty,$ the last formula turns into 
 $$\langle  \partial_ t u, u-\xi\rangle  =\frac{1}{2}
 \frac{d}{dt}  \int_{\Omega}  \vert u-\xi\vert^2\: dx  .$$ 
This is the common term for the standard notion of variational solution. It is noteworthy, however,   that the truncation operation here is an important ingredient to get uniqueness (see the uniqueness proof and the remark below).

 \bigskip 
The considerations above  bring on  the following definition of variational  solutions :
 % Definition 
\begin{defi}\label{defweaksol} 
	Let $\displaystyle f\in L^{1}( \Omega_T)$ and $u_0\in   Lip_1.$  A variational  solution of  (E) is a function   $u\in L^{\infty}(0,T, C_0(\Omega))$ such that $u(t) \in   Lip_1$ for a.e. $t\in [0,T]$,  and for every $\xi \in Lip_1$  and every $k>0$,
	\begin{equation}\label{varform}
	\frac{d}{dt}  \int_{\Omega}  \int_0^{u(t)} T_k (s -\xi) ds  { d}x  - \int_{\Omega} \gamma(u)\: \mathcal{H} ( \partial_x K   \ast u )  \partial_xT_k (u-\xi)  dx 
	\le  \int_{\Omega} f\:  T_k(u-\xi) \:  dx 
	\end{equation} 
	in $\mathcal D'([0,T)).$  
\end{defi}

In other words,  $u\in L^{\infty}((0,T), C^0(\Omega))$ such that $u(t) \in  Lip_1$ for a.e. $t\in [0,T]$ is a variational  solution of (E) if for every $\xi  \in Lip_1$ and every $\sigma \in C^1([0, T), \RR_+)$ one has
$$- \int_0^T\!\!\!  \int_{\Omega} \dot{\sigma}(t) \int_0^{u(t)} T_k (s -\xi) ds   dx  dt 
- \int_0^T\!\!\!  \int_{\Omega}   {\sigma}(t) \gamma(u(t))\: \mathcal{H} ( \partial_x K   \ast  )  \partial_xT_k (u(t)-\xi)  dx   dt  $$  $$
 \leq  \sigma(0) \int_{\Omega} \int_0^{u_0} T_k (s -\xi) ds  d x   + \int_0^T\!\!\!  \int_{\Omega} \sigma(t) f (t)\: T_k(u(t)-\xi)    d t d x.$$

%Remark  
\begin{rem}\label{rsol}
\begin{enumerate}
	\item  It is achievable    to 	define a solution as  a function  $u\in L^{\infty}((0,T), C^0(\Omega))$,  with $\partial_t u \in L^1(0,T; Lip^*_0)$ and $u(t) \in  Lip_1$ for a.e. $t\in [0,T]$,  $u(0)=u_0$ 
	and 
		\begin{equation}\label{wf01}
	\begin{array}{c} 	 \langle  \partial_ t u, T_k(u-\xi) \rangle   - \int_{\Omega} \gamma(u)\: \mathcal{H} ( \partial_x K   \ast u )  \partial_xT_k (u-\xi)  dx \\  
		\le  \int_{\Omega} f\:  T_k(u-\xi) \:  dx   ,\quad \hbox{ for a.e. }t\in (0,T),
		\end{array}
	\end{equation} 
	where $\langle .,.\rangle$ denote sthe duality brackect in $Lip^*_0.$ 
However, one can prove that  if $u$ satisfies (\ref{wf01})  it is also a variational solution  in the sense of Definition \ref{defweaksol}.  Indeed, if  $\partial_t u \in L^1(0,T; Lip^*_0)$ one can prove rigorously   that  (\ref{wf01})  yields 
	$$ \int _0^T \langle  \partial_ t u(t), T_k(u(t)-\xi) \rangle \sigma(t)\:  dt  =
	- \int_0^T \int_{\Omega} \dot{\sigma}(t) \int_0^{u(t)} T_k (s  -\xi) ds   dx  dt $$ 
	$$-    \sigma(0) \int_{\Omega} \int_0^{u_0} T_k (s -\xi) ds  d x  ,$$ 
	for any $\xi \in Lip_1$ and $\sigma \in C^1([0, T), \RR).$

\item Notice that some similar notion of solution have been used in \cite{AgCaIg} for a different problem using the  so called W1-JKO scheme, where     W1 is related to the Wasserstein distance $W_1.$  
\end{enumerate} 
\end{rem}

 % Theorem 
\begin{thm}\label{theo0}
Let $\displaystyle f\in L^{1}( \Omega_T)$ and $u_0\in   Lip_1.$  The problem (E)  has a unique variational solution $u$.
 \end{thm}

\bigskip
To prove this theorem, we see that  Lemma \ref{ldefsol} implies that the operator $u\in Lip_1\to -\nabla \cdot (m\nabla u)$,   with non-negative  $m$ satisfying $m(\vert \nabla u\vert -1)=0,$  may be represented in 	$L^2(\Omega),$ by  the sub-differential operator $\partial \ind_{{Lip_1}}$ of the indicator function $  \ind_{{Lip_1}}\: :\:  L^2(\Omega) \to [0,\infty],$  
	$$\displaystyle  \ind_{{Lip_1}} (z)= \left\{ \begin{array}{ll} 
	0\quad &\hbox{ if  } z\in {Lip_1} \\  \\
	\infty &  \hbox{ otherwise .}\end{array}
	\right.  $$
In particular, this implies that the equation (E) is formally of the type 
\begin{equation}\label{Etype}
\displaystyle \frac{du}{dt}+\partial \ind_{{Lip_1}} u_{}\ni \T(u) +f\quad  \hbox{ in }(0,T),
\end{equation} 
where 
  $\T\: :\:  {Lip_1}\subset L^2(\Omega) \to L^2(\Omega)$ is given by  
$$\T(u)= -  \partial_x \Big( \gamma(u) \mathcal{H}(  K \star   \partial_x u )\Big),\quad \hbox{ for any }u\in {Lip_1}.    $$

Recall that the case where $\mathcal{H}\equiv 0,$ the phenomena  corresponds  simply to the sandpile 
problem where the dynamic is completely governed by the following nonlinear evolution equation : 
\begin{equation}\label{dyn1}
\left\{ \begin{array}{ll} 
\displaystyle \frac{du}{dt}+\partial \ind_{{Lip_1}} u\ni f\quad &\hbox{ in }(0,T)\\  \\
\displaystyle u(0)=u_0, \end{array}
\right.
\end{equation}
in $L^2(\Omega).$

\bigskip
For the proof of Theorem  \ref{theo0}, we begin with the following results concerning (\ref{dyn1}) which will be  useful.
%Proposition
 \begin{pro} \label{pV0}
 For any $f\in L^2(\Omega_T)$ and $u_0\in Lip_1,$  there exists a unique solution of the problem (\ref{dyn1}), in the sense that  $u\in W^{1,\infty}(0,T; L^2(\Omega)),$ 
 $u(0)=u_0$ and $$ \displaystyle f(t)- \frac{du(t)}{dt}\in \partial \ind_{{Lip_1}} u(t) \quad \hbox{ for a.e.  }t\in (0,T).$$ 
 Moreover, we have 
 \begin{enumerate} 
 \item \label{P1} $u\in L^\infty(0,T;\C_0^{0,\alpha}(\Omega))\cap L^p(0,T;W^{1,p}(\Omega)),$ for any $0\leq \alpha<1$ and $1\leq p<\infty,$   and $u(t)\in {Lip_1}$ for a.e.  $t\in [0,T).$  
 \item \label{P2} $\partial_t u\in L^1(0,T;Lip^*_0)$ and we have 
\begin{equation}\label{Lip'}
 \Vert \partial_t u\Vert_{L^1(0,T;Lip^*_0)}  \leq  2\: \delta_\Omega \Vert f\Vert_{L^1(0,T;Lip^*_0)}     + \frac{1}{2}  \int u_{0}^2.
 \end{equation}
  \end{enumerate}
  \end{pro}
 %%%%%%%
\noindent {\textbf{Proof:}}  The existence of a solution  $u\in W^{1,\infty}(0,T; L^2(\Omega))$ follows by standard theory of evolution problems governed by sub-differential operator (cf. \cite{Br}).  By definition of the solution, we know  that $u(t)\in {Lip_1}$   and $\vert u(t) \vert \leq \delta_\Omega$ in $\Omega$, for any $t\in [0,T).$    Using the fact that  $ {Lip_1}$ is compactly injected in $\C_0^{0,\alpha}(\Omega),$ we deduce that  $u\in L^\infty(0,T;\C_0^{0,\alpha}(\Omega))$.   Thus (\ref{P1}).  Let us  prove (\ref{P2}). For any $\xi\in Lip_1,$ we see that testing with $-\xi$ and letting $k\to\infty,$ we have 
$$ \int_\Omega \left( f(t)-\partial_t u(t) \right) (u(t)+\xi)\: dx \geq 0,\quad \hbox{ for any } t\in [0,T).$$ 
This implies that 
$$\int_\Omega  \partial_t u(t) \: \xi \: dx \leq \int_\Omega  f(t)(u(t)+\xi)\: dx -\frac{1}{2}\frac{d}{dt}   \int_\Omega u(t)^2\: dx .$$
Integrating over $(0,T),$ we get  
\begin{eqnarray*} 	\int_0^T\!\! \int_\Omega \partial_t u (t)\:  \xi   \: dtdx  &\leq&   \int_0^T\!\! \int_\Omega   f\:( \xi   +  u)    + \frac{1}{2}  \int u_{0}^2(t) \: dx\\  \\
	&\leq& 2\: \delta_\Omega \Vert f\Vert_{L^1(0,T;Lip^*_0)}     + \frac{1}{2}  \int u_{0}^2\: dx.
	\end {eqnarray*}
Since $\xi$ is arbitrary in $Lip_1,$   we deduce    (\ref{Lip'}).   

\qed

Now, coming back to the problem (\ref{Etype}),  thanks to the assumptions on $\mathcal{H},$  $K$ and $\gamma,$ the operator $\T$ is well defined, and for any   $z\in L^2(0,T;W^{1,2}(\Omega)),$  we have 
$\T (z)\in L^2(Q).$  So,  given    $u_0\in {Lip_1},$  thanks to Proposition \ref{pV0}, the sequence $(u_n)_{n\in\NN} $ given  by 
\begin{equation}\label{dyn1n}
 \displaystyle \frac{du_{n+1}}{dt}+\partial \ind_{{Lip_1}} u_{n+1}\ni \T(u_n) +f\quad  \hbox{ in }(0,T),\:   \hbox{ for }n=0,1,2 ...
 \end{equation}
 is well defined in $W^{1,\infty}(0,T; L^2(\Omega)) \cap L^p(0,T;W^{1,p}(\Omega)),$ for any   $1\leq p<\infty.$    Moreover, we have

%Lemma
\begin{lem}\label{ln} 
\begin{enumerate} 
\item $u_n$ is a bounded sequence in $L^\infty(0,T;\C_0^{0,\alpha}(\Omega)),$ for  $0\leq \alpha<1.$
%\item $u_n$ is a bounded sequence in $L^p(0,T;W^{1,p}_0( \Omega)),$ for anu  $p\geq 1.$
\item $\partial_t u_n$ is a bounded sequence in $L^1(0,T;Lip^*_0).$ 
 \end{enumerate} 
 \end{lem}
 %%%%%%%
\noindent {\textbf{Proof:}} \begin{enumerate} 
\item Thanks  to Proposition \ref{pV0},  we know that   $u_n\in L^\infty(0,T;\C_0^{0,\alpha}(\Omega)),$ $0\leq \alpha<1,$  and $u_n(t)\in {Lip_1}$ for a.e.  $t\in [0,T).$ This implies that  $u_n$ is bounded in $L^\infty(0,T;\C_0^{0,\alpha}(\Omega)).$

\item  Thanks again to Proposition \ref{pV0},  we have 
\begin{eqnarray*} 
 \Vert \partial_t u_{n+1}\Vert_{L^1(0,T;Lip^*_0)}   
&\leq &  2  \: \delta_\Omega \:\Vert \T (u_n)\Vert_{L^1(0,T;Lip^*_0)}  +2 \: \delta_\Omega  \:\Vert f \Vert_{L^1(0,T;Lip^*_0)}   + \frac{1}{2}  \int u_{0}^2\\  \\ 
&\leq& 2\: \delta_\Omega   \:\left \Vert \partial_x \Big( \gamma(u_n) \mathcal{H}(  K \star   \partial_x u_n )\Big) \right \Vert_{L^1(0,T;Lip^*_0)} +2 \: \delta_\Omega \:\Vert f \Vert_{L^1(0,T;Lip^*_0)} \\  \\ 
&  & + \frac{1}{2}  \int u_{0}^ 2\\  \\ 
&\leq& 2  \delta_\Omega   \:   \left \Vert    \gamma(u_n) \mathcal{H}(  K \star   \partial_x u_n )  \right \Vert_{L^1(Q_T)} +2  \delta_\Omega   \:\Vert f \Vert_{L^1(0,T;Lip^*_0)}  + \frac{1}{2}  \int u_{0}^ 2 .
 \end{eqnarray*} 
Using the fact that $\gamma$ and $\mathcal H$ are Lipschitz continuous and that $ u_n(t)\in {Lip_1},$ we deduce that there exists $C$ (independent of $n$) such that 
$$ \Vert \partial_t u_{n+1}\Vert_{L^1(0,T;Lip^*_0)}    \leq C. $$
Thus the result of the lemma.   
  \end{enumerate} 
\qed

\bigskip 
\noindent {\bf{Proof of Theorem \ref{theo0} :}} 
%%%%%%%%%%%%%%%%%
 \noindent \underline{\bf{Existence :}}  First assume that $f\in L^2(\Omega_T).$ Let us consider the sequence $(u_n)_{n\in \NN}$ as given by Lemma \ref{ln}.  Since the embedding $\C_0^{0,\alpha}(\Omega)$ into $\C_0(\Omega)$ is compact and the   embedding of $\C_0(\Omega)$ into $Lip^*_0$ is continuous, by using  Lemma 9 of \cite{simon}, we can conclude that, by taking a subsequence if necessary, $u_n$ converges to $u$ in $L^1(0,T; \C_0(\Omega)),$ and we have  $u(t)\in {Lip_1},$ for a.e. $t\in [0,T).$ 
Since, for a.e. $t\in [0,T)$ and any $\xi\in Lip_1,$   $u_n(t)-T_k(u_n(t) -\xi)\in Lip_1,$  (\ref{dyn1n}) implies that 
$$\int_\Omega \frac{\partial u_{n+1}(t)}{\partial t}\: T_k(u_{n+1} (t) -\xi)  
-    \int_\Omega      \gamma(u_n) \mathcal{H}(  K \star   \partial_x u_n )\: \partial_x T_k(u_{n+1} (t) -\xi)  $$   $$\leq \int f\:  T_k(u_{n+1} (t) -\xi) ,$$ 
so that 
$$ \frac{d}{dt}  \int_{\Omega}  \int_0^{u_{n+1}(t)} T_k (s -\xi) ds  dx  -  \int_{\Omega} \gamma(u_{n})\: \mathcal{H} ( \partial_x K   \ast u_{n} )  \partial_x T_k(u_{n+1}-\xi)  dx 
  $$  $$\le  \int_{\Omega} f\:  T_k(u_{n+1}-\xi) \:  dx $$
in $\mathcal D'([0,T).$  
Then letting $n\to \infty,$ and using the convergence of $u_n$ in $L^1(0,T; \C_0(\Omega))$ and Lebesgue dominated convergence theorem, we obtain (\ref{varform}).   Now for $\displaystyle f\in L^{1}( \Omega_T)$ we consider   $\displaystyle f_m\in L^{2}( \Omega_T)$ such that $ \displaystyle f_m \to f $ in $\displaystyle L^{1}( \Omega_T)$  and the sequence $(u_m)_{m\in \NN}$  given  by 
\begin{equation}\label{dyn1m}
\displaystyle \frac{du_{m}}{dt}+\partial \ind_{{Lip_1}} u_{m}  +\partial_x\big( \gamma(u_m) \mathcal{H}(  K \star   \partial_x u_m ) \big)\ni f_m\quad  \hbox{ in }(0,T),\:   \hbox{ for }m=0,1,2 ...
\end{equation}
in the following sense
\begin{equation}\label{varformm}
\frac{d}{dt}  \int_{\Omega}  \int_0^{u_m(t)} T_k (s -\xi) ds  dx  - \int_{\Omega} \gamma(u_m)\: \mathcal{H} ( \partial_x K   \star u_m )  \partial_xT_k (u_m-\xi)  dx 
\le  \int_{\Omega} f_m\:  T_k(u_m-\xi) \:  dx 
\end{equation} 
in $\mathcal D'([0,T)$  for any $\xi \in Lip_1$  and $k>0$. Similarly as in lemme \ref{ln}, we have
\begin{eqnarray*} 
	\Vert \partial_t u_{m}\Vert_{L^1(0,T;Lip^*_0)}   \leq2 \delta_\Omega\: \left \Vert    \gamma(u_m) \mathcal{H}(  K \star   \partial_x u_m )  \right \Vert_{L^1(\Omega_T)} +2\delta_\Omega  \:\Vert f_m \Vert_{L^1(0,T;Lip^*_0)}  + \frac{1}{2}  \int u_{0}^ 2 dx 
\end{eqnarray*} 
and $u_m \to u $ in $L^1(0,T; \C_0(\Omega))$. Letting $m\to \infty,$ and using the
dominated convergence theorem,  the proof of the existence  is finished.\\
  \\ \\
  %%%%%%%%%%%%%%%%%
 \noindent \underline{\bf{Uniqueness :}}   Now, to prove  the uniqueness, let $\displaystyle u_{1}$ and $u_{2}$
 be two solutions of $(E)$ in the sense of (\ref{varform}). We have 
  \begin{equation}
 \frac{d}{dt}  \int_{\Omega}  \int_0^{u_i(t)} T_n (s -\xi) ds  dx  - \int_{\Omega} \gamma(u_i)\: \mathcal{H} ( \partial_x K   \ast u_i )  \partial_xT_n (u_i-\xi)  dx 
 \le  \int_{\Omega} f\:  T_n(u_i-\xi) \:  dx 
 \end{equation} 
 To double variables, we consider $u_1=u_1(t)$ and $u_2=u_2(s),$ for any $s,t\in [0,T).$   Using the fact that $u_1=u_1(t)$  is a solution  and setting $\xi =u_2(s) $  which is considered constant with respect to $t$      have 
 $$  \frac{d}{dt}  \int_{\Omega}  \int_0^{u_1(t)} T_n (r -u_2(s)) dr  dx  - \int_{\Omega} \gamma(u_1(t))\: \mathcal{H} ( \partial_x K   \ast u_1(t))  \partial_xT_n (u_1(t)-u_2(s))  dx $$ $$
 \le  \int_{\Omega} f(t)\:  T_n(u_1(t)-u_2(s)) \:  dx .$$
In the same way, taking    $u_2=u_2(s)$  is a solution  and setting $\xi =u_1(t)$,  we have 
  $$  \frac{d}{ds}  \int_{\Omega}  \int_0^{u_2(s)} T_n (r -u_1(t)) dr dx  - \int_{\Omega} \gamma(u_2(s))\: \mathcal{H} ( \partial_x K   \ast u_2(s))  \partial_xT_n (u_2(s)-u_1(t))  dx $$ $$
 \le  \int_{\Omega} f(s)\:  T_n(u_2(s)-u_1(t)) \:  dx \Big). $$
  Dividing by $n,$ and adding the two equations, we obtain 
  $$ \frac{1}{n} \frac{d}{dt}  \int_{\Omega}  \int_0^{u_1(t)} T_n (r -u_2(s)) dr  dx  +  
  \frac{1}{n}    \frac{d}{ds}  \int_{\Omega}  \int_0^{u_2(s)} T_n(r -u_1(t)) dr dx   $$ 
  $$ \leq  \frac{1}{n}  \int_{\Omega} \Big\{ \gamma(u_1(t))\: \mathcal{H} ( \partial_x K   \ast u_1(t))       -\gamma(u_2(s))\: \mathcal{H} ( \partial_x K   \ast u_2(s)) \Big\} \partial_xT_n (u_1(t)-u_2(s))\:  dx $$
  $$+  \frac{1}{n}   \int_{\Omega} ( f(t) -f(s)) \:  T_n(u_1(t)-u_2(s)) \:  dx   $$
   Let us re-write the second equation  
  $$ I_n:=  \frac{1}{n}\int_{\Omega} \Big\{ \gamma(u_1(t))\: \mathcal{H} ( \partial_x K   \ast u_1(t))       -\gamma(u_2(s))\: \mathcal{H} ( \partial_x K   \ast u_2(s)) \Big\} \partial_xT_n (u_1(t)-u_2(s))\:  dx  $$
  as \begin{equation}\label{Ex445} I_n  =:I_n^1 +I_n^2,    \end{equation} 
  with 
  $$I_n^1:=  \frac{1}{n}\:   \int_{\Omega}  (\gamma(u_1(t))- \gamma( u_2(s))) \partial_x\: T_n ( u_1(t)-  u_2 (s) ) \mathcal{H} ( \partial_x K   \ast u_1(t) )  \:dx  $$
  and 
  $$ I_n^2:=   \frac{1}{n}\:   \int_{\Omega}   \gamma(u_2(s))\Big(\mathcal{H} ( \partial_x K   \ast u_1(t) )- \mathcal{H} ( \partial_x K   \ast u_2(s) )\Big)  \partial_x\: T_n ( u_1(t)-  u_2(s)  ) \:dx .$$ 
  Recall that $u_i(t)\in Lip_1,$ for any $t\in [0,T),$ ]$\gamma$ and $\mathcal H $ are  Lipschitz continuous.  So, there exists a constant $c^\star>0$ (independent of $n$), such that  
  \begin{equation}\label{Ex445}
  I_n^1  \leq  c^\star\:  \int_{\Omega}  |u_1(t)-  u_2(s)| \:dx .
  \end{equation}
  Integrating  by parts in $I_n^2,$  we obtain
  \begin{equation}\label{Ex5}
  \begin{array}{l l }
  \displaystyle       I_n^2 
  & \displaystyle=   - \frac{1}{n}\:  \int_{\Omega}   \gamma^\prime(u_2(s)) \partial_x u_2(s)\Big(\mathcal{H} ( \partial_x K   \ast u_1(t) )- \mathcal{H} ( \partial_x K   \ast u_2(s) )\Big)  \: T_n ( {u_1(t)-  u_2(s)} ) \:dx   \\    \\  
  & \displaystyle  - \frac{1}{n}\:   \int_{\Omega}   \Big\{   \gamma(u_2(s)) \Big(  \big( \partial_{x^2} K   \ast u_1(t)\big) \mathcal{H}^\prime ( \partial_x K   \ast u_1(t) )  \\  
  & 
  \hspace*{4cm}
  -     \big(\partial_{x^2} K   \ast u_2(s)\big) \mathcal{H}^\prime ( \partial_x K   \ast u_2(s) )\Big)  \: T_n ( {u_1(t)-  u_2(s)} )\Big\} \:dx . 
  \end{array} 
  \end{equation}
 The first term   of $I_n^2$ satisfies 
  \begin{equation}\label{Ex556}
  \begin{array}{l}  \displaystyle \frac{1}{n}\:   \Big| \int_{\Omega}    \gamma^\prime(u_2(s)) \partial_x u_2(s)\Big(\mathcal{H} ( \partial_x K   \ast u_1(t) )- \mathcal{H} ( \partial_x K   \ast u_2(s) )\Big)  \: T_n ( {u_1(t)-  u_2(s)} ) \:dx   \Big| \\  \\ 
  \hspace*{3cm} \displaystyle
  \leq  c^{\star\star}\: \int_{\Omega}  |u_1(t)-  u_2(s)| \:dx .
  \end{array}
  \end{equation}
  As to the second term  that we denote here by $I_n^{2'}$ 
  $$\displaystyle  I_n^{2'} :=\frac{1}{n}\: \int_{\Omega}    \gamma(u_2(s)) \Big\{  \big( \partial_{x^2} K   \ast u_1(t)\big) \mathcal{H}^\prime ( \partial_x K   \ast u_1(t) )- $$  $$\hspace*{5cm}\big(\partial_{x^2} K   \ast u_2(s)\big) \mathcal{H}^\prime ( \partial_x K   \ast u_2(s) )\Big\}  \: T_n ( u_1(t)-  u_2(s) ) \:dx   ,$$
 we have  
  \begin{equation}\label{Ex66}
  \begin{array}{l l }
  \displaystyle  |I_n^{2'}|  \leq  c^{\star\star \star}\:  \int_{\Omega}  |u_1(t)-  u_2(s)| \:dx  .
  \end{array} 
  \end{equation}
  From (\ref{Ex445}), (\ref{Ex556}) and (\ref{Ex66}), we obtain 
  $$ \vert I_n\vert \leq       \: C\:  \int_{\Omega}  |u_1(t)-  u_2(s)| \:dx  ,$$
  so that, for any $n>0,$ we have   
    $$  \frac{d}{dt}  \int_{\Omega}  \int_0^{u_1(t)} T_n (r -u_2(s)) dr  dx  +  
  \frac{d}{ds}  \int_{\Omega}  \int_0^{u_2(s)} T_n(r -u_1(t)) dr dx   $$ 
  $$ \leq  \: C\:  \int_{\Omega}  |u_1(t)-  u_2(s)| \:dx  +  \int_{\Omega} ( f(t) -f(s)) \:  T_n(u_1(t)-u_2(s)) \:  dx  . $$
  Letting $n\to 0,$  we get 
   $$  \frac{d}{dt}  \int_{\Omega}  \int_0^{u_1(t)} \mbox{sign}_0 (r -u_2(s)) dr  dx  +  
  \frac{d}{ds}  \int_{\Omega}  \int_0^{u_2(s)} \mbox{sign}_0(r -u_1(t)) dr dx   $$ 
  $$ \leq  \: C\:  \int_{\Omega}  |u_1(t)-  u_2(s)| \:dx  +  \int_{\Omega} ( f(t) -f(s)) \:  \mbox{sign}_0(u_1(t)-u_2(s)) \:  dx  . $$ 
  Thus   $$
  \displaystyle  \frac{d}{dt} \int_{\Omega}  | u_1(t)-u_2(s)| \:dx  +  \frac{d}{ds} \int_{\Omega}  | u_1(t)-u_2(s)| \:dx  \leq   
  \displaystyle   \nu  \: C\: \int_{\Omega}  |u_1(t)-  u_2(s)| \:dx  $$  $$+  \int_{\Omega} \vert  f(t) -f(s)) \vert \:  dx.
  $$
  Now, de-doubling variables  $t$ and $s,$  we get   
 \begin{equation}\label{inter2}\frac{d}{dt}  \int_{\Omega } \vert u_1(t)-u_2(t)\vert   dx
 \le C   \int_{\Omega}  \vert  u_1-u_2\vert       dx     \quad \hbox{ in }\mathcal{D}'(0,T) 
 \end{equation} 
 and the uniqueness follows by  Gronwall Lemma.
 
\qed

 \begin{rem}
 	See in the proof of   Theorem \ref{theo0}, that   it is possible to prove the result of existence  of a variational solution for any $f\in L^1(0,T,\left( \C_0^{0,\alpha}(\Omega)\right)^*),$ with $0\leq \alpha<1,$ where  $\left( \C_0^{0,\alpha}(\Omega)\right)^*$ denotes the topological dual space of $\C_0^{0,\alpha}(\Omega).$ More precisely, 
for any $f\in L^1(0,T, \left( \C_0^{0,\alpha}(\Omega)\right)^* )$ and $u_0\in Lip_1,$ the problem   (E)  has a variational solution in the sense that   $u\in L^{\infty}(0,T, C_0(\Omega))$,  $u(t) \in   Lip_1$ for a.e. $t\in [0,T]$,  and for every $\xi \in Lip_1$  and every $k>0$,
$$\frac{d}{dt}  \int_{\Omega}  \int_0^{u(t)} T_k (s -\xi) ds  { d}x  - \int_{\Omega} \gamma(u(t))\: \mathcal{H} ( \partial_x K   \ast u(t) )  \partial_xT_k (u(t)-\xi)  dx  $$  $$
\le \langle f(t), T_k(u(t)-\xi) \rangle, \quad \hbox{ in } \mathcal D'([0,T)).$$
This allows in particular to consider the situations where we have some singular source terms of the type 
$$f(t)=\sum_n \left( \delta_{x_n} - \delta_{y_n}   \right),  $$
where $x_n$ and $y_n$ are sequences in $\RR^d,$ satisfying 
$$\sum_n \vert x_n-y_n\vert ^\alpha <\infty.   $$   
However, the uniqueness is not clear if one weaken the assumption $f\in L^1(\Omega_T).$ 
  \end{rem}

  \section*{Acknowledgements} 
 \noindent This work was performed under the research project PPR CNRST : Mod\`eles Math\'ematiques appliqu\'es \`a l'environnement, \`a
l'imagerie m\'edicale et aux Biosyst\`emes (Essaouira-Morocco).

\end{document}